\documentclass[a4paper,10pt]{article}
\usepackage{kubik}

\def\rem#1{\par\medskip\noindent\refstepcounter{thm}\hbox{\bf \arabic{section}.\arabic{thm}. #1.}
\sf \ %\ignorespaces
}

\newenvironment{rdef}{\rem}

\begin{document}
 
\title{Alexandrov spaces with\\ 
maximal number of extremal points}
\author{Nina Lebedeva}\maketitle

\begin{abstract}
We show that any $n$-dimensional
nonnegatively curved Alexandrov space with the maximal possible number of extremal points
is isometric to  a quotient space of $\RR^n$   by an action of a crystallographic group.
We  describe all such actions. 
\end{abstract}

\footnotetext[1]{Partially supported by RFBR grant 11-01-00302a}

\section{Introduction}

If the space of directions at a  point $p$ in an Alexandrov space 
has diameter $\le \tfrac\pi2$, this point
is called \emph{extremal}.
 Equivalently, the one-point set $\{p\}$ is an \emph{extremal set} as  defined by Perelman and Petrunin in \cite{perelman-petrunin:extremal}.
Yet equivalently, 
$p$ is a \emph{critical point} of every distance function.

It has been proven by Perelman (\cite{Per94}) that every
 $n$-dimensional Alexandrov space with nonnegative
curvature has at most $2^n$ extremal points.
For completeness, we present this proof in Subsection~\ref{subsec:2^n}.
This proof is a slight modification of a proof of the following problem in discrete geometry:

\begin{thm}{Problem}
Assume $x_1,x_2,\dots,x_m$ is a collection of points in the $n$-dimensional Euclidean space such that $\angle x_ix_jx_k\le \tfrac\pi2$ for any distinct $i$, $j$ and $k$. 
Show that $m\le 2^n$ and moreover, if $m= 2^n$ then the $x_i$ form 
the set of vertices of a right parallelepiped.
\end{thm}

This problem 
%had been
 posted by Erd\H os in \cite{erdos} was solved by Danzer and Gr\"unbaum in \cite{danzer-gruenbaum}.

\medskip

In this paper we study 
nonnegatively curved $n$-dimensional Alexandrov spaces with $2^n$ extremal points,
we call such spaces  \emph{$n$-boxes}.

Classification of  $n$-boxes is a folklore problem.
Clearly, right parallelepipeds are boxes. It was
% In a private conversation around 1993, G.~Perelman
 suggested (\cite{Per93}) that these
 might be the only examples.
Soon it was noticed (\cite{ Kap}) that %this condition also holds for 
the boundary of the  $3$-dimensional Euclidean  tetrahedron
whose opposite edges are equal
(or equivalently, whose four faces are congruent  triangles)
 is also a $2$-box.
Latter, it was conjectured (\cite{Pet}) that all $n$-boxes have to be isometric to a
quotient 
 of a flat torus by an action of a group of isometries which is isomorphic to a product of $\ZZ_2$-groups.
However, it turns  out
  that  not all $n$-boxes can be obtained this way.
  The first counterexample arrises in dimension 3, this is a space $\text{\mancube}_2'$ constructed below in this section.

Our main results are Theorem~\ref{main1} and 
Theorem~\ref{main2}.

\begin{thm}{Theorem}\label{main1}
For any $n$-box
there exists a group $\Gamma$
and a discrete cocompact isometric action
$\Gamma\acts\RR^n$
 such that this $n$-box 
 is isometric to the quotient space $\RR^n/\Gamma$.
\end{thm}

Theorem~\ref{main2} below describes  all possible actions on $\RR^n$ which 
produce $n$-boxes.
Proposition~\ref{prop:affine+extreme} 
implies that it is sufficient to describe the actions $\Gamma\acts\RR^n$ up to affine conjugation.
We need firstly the following
\begin{rdef}{Definition}
Let $\Gamma\acts\RR^n$ be a group action.
We call
a point $x\in \R^n$ a \emph{singular point} (for the action $\Gamma\acts\RR^n$)
if it is a unique fixed point for some subgroup of $\Gamma$.
\end{rdef}

\begin{thm}{Proposition}\label{prop:affine+extreme}
For  a discrete action $\Gamma\acts\RR^n$ by isometries 
the quotient space 
$\mathcal A=\RR^n/\Gamma$ is an Alexandrov space
of nonnegative curvature.
 Moreover, a  point $e\in \mathcal A$ is extremal if and only if it is an 
image of a singular point.
%an isolated fixed point in $\RR^n$ of some subgroup of $\Gamma$. 

In particular
if two such actions are affine conjugate, then
 the number of extremal points in corresponding quotient spaces are equal.
\end{thm}
\parit{Proof.}
Clearly $\mathcal A$ is a polyhedral space and any angle around any face of codimension 2
in $\mathcal A$ is less or equal $\pi/2$, hence
$\mathcal A$ is an Alexandrov space
of nonnegative curvature.
Obviously, the image of any singular point is a vertex in the polyhedron $\mathcal A$,
and any pre-image of any vertex in $\mathcal A$ is a singular point.
It remains to note that for any discrete isometric action $G\acts S^{n-1}$
a condition $\diam S^{n-1}/G<\pi$ implies that  $\diam S^{n-1}/G\le\pi/2$, hence
 all vertices in $\mathcal A$ are extremal points.
\qeds

Recall that the Coxeter group associated
with an
  $n$-polyhedron is the group generated by reflections in its faces;
such a group is defined together with an action on $\RR^n$.
Let us denote by $\bigoplus^n\acts\RR^n$ the action of the Coxeter group
$\bigoplus^n$  of the unit cube.

\begin{thm}{Theorem}\label{main2}
Let $\Gamma\acts\RR^n$ be a subaction of $\bigoplus^n\acts\RR^n$ such that
 % the stabilizer of
  any vertex $e$ of the unit cube is an isolated fixed point for some subgroup of $\Gamma$.
Then $\RR^n/\Gamma$ is an $n$-box.

Moreover, all  $n$-boxes
arise from such 
actions
$\Gamma\acts\RR^n$ 
  or their affine conjugate actions.
\end{thm}

It is immediate from the theorem that $[\bigoplus^n:\Gamma]=2^k$ for some $k\in\{0,\dots,n-1\}$.
Note that Theorem~\ref{main2} makes possible to list all  group actions
which produce $n$-boxes.
Let us fix a set $S$ of faces
of the $n$-cube $Q$ with the following property.
Any vertex of $Q$ is the only intersection of all faces in $S$ containing this vertex.
Then the group action generated by reflections in the elements of 
$S$  gives an example of 
an action $\Gamma\acts\RR^n$ 
described in Theorem~\ref{main2} 
and any such an action can be obtained in this way.
It remains to find all isometric actions that are affine conjugate to the constructed action.
This is equivalent to finding all parallel metrics invariant w.r.t. reflections in the elements of $S$.
Therefore, 
one can think of
 any $n$-box
 as a space
   glued from $2^k$ copies of the cube equipped with a parallel metric $g$ which is invariant under  all reflections in  the elements of $S$. 

Let us use our construction to classify $n$-boxes in low dimensions.
\begin{itemize}
\item For $n=1$ 
there exists only one, up to affine conjugation group action, 
which produces $1$-box. The corresponding quotient space
 $\II=[0,1]$ cares a $1$-parametric family of metrics.
\item For $n=2$ there exist two spaces (up to a choice of a parallel metric): the square $\square=\II\times\II$ and the double square $\square_2$. 
The square $\square$ admits a $2$-parametric family of metrics;
this family gives rise to 
 all possible rectangles. 
The double square $\square_2$ admits a 3-parametric family of metrics;
this family gives rise to surfaces of 3-simplexes whose opposite (non-intersecting) edges are equal.
Such simplexes are sometimes called \emph{disphenoids}.
\item For $n=3$ there are five $3$-boxes (up to a choice of a parallel metric): 
the cube $\text{\mancube}=\II\times\II\times\II$; 
the double cube $\text{\mancube}_2$ (obtained by gluing two copies of the cube along their common boundary);
doubling of the cube in the 5 faces $\text{\mancube}_2'$
 (obtained by gluing two copies of the cube along 5 faces of their common boundary);
the product $\text{\mancube}_2''=\II\times \square_2$
and the quotient $\text{\mancube}_4$ of the standard torus by the central symmetry.
The dimensions of the space of metrics are respectfully
$3$, $3$, $3$, $4$ and $6$.
\end{itemize}

%A big part of my proof works for this conjecture, but there is a apparently small gap which I can not pass,
%see my question on mathoverflow \cite{mathoverflow}.

\parbf{Structure of the paper.}
In Subsection~\ref{subsec:2^n}
for the sake of completeness we
reproduce
 the proof
that $n$-dimensional Alexandrov space with nonnegative
curvature has at most $2^n$ extremal points.

The proof of
Theorem~\ref{main1} (Sections 2,3,4) is organized in two steps.

In Sections 2,3 we show that an $n$-box $\mathcal A$ has to be a polyhedral space (Theorem~\ref{thm:polih}).
According to Proposition~\ref{prop:poly-char},
it is sufficient to show that each point $p\in \mathcal A$ has a conic neighborhood (see Definition~\ref{con.negh}).
This is proved in Key lemma~\ref{lem:cone}.

In Section 4, we show that $\mathcal A$ is a flat orbifold.
By Proposition~\ref{prop:orbi-char}, it is sufficient to show 
that an
angle around any face of codimension $2$ in $\mathcal A$ has to be $\pi$, or $2\pi$.
This is proved in Theorem~\ref{cone:::pi}.

In sections 5,6,7
we prove  Theorem~\ref{main2}.

\parbf{Acknowledgements.}
I thank Anton Petrunin for
bringing this problem to my attention  and for useful discussions.

\subsection{The upper estimate for the number of extremal points}
\label{subsec:2^n}

In this subsection we give the proof (due to Erd\H{o}s, Danzer, Gr\"unbaun, Perelman)
of Theorem~\ref{thm:2^n}.
We also introduce notation which are used further.

\begin{thm}{Theorem}\label{thm:2^n}
The number of extremal points of an $n$-dimensional nonnegatively curved Alexandrov space is at most $2^n$.
\end{thm}
\begin{rdef}{Notation}
Denote by $ A$  an $n$-dimensional 
nonnegatively curved Alexandrov space.

We label the extremal points in $ A$ by $e_1,e_2,\dots,e_m$.

For a triangle $abc$ in    $ A$
we denote by
$\~a\~b\~c$ a comparison triangle in $\R^2$ (i.e. a triangle with the same lengths of sides).
We denote 
 by $\~\angle abc$ the angle at $\~b$ of the triangle
$\~a\~b\~c$.

For a point $ a\in A$
we denote by $\Sigma_a$ the unite tangent space at $a$.

\end{rdef}

\parit{Proof of  {Theorem}~\ref{thm:2^n}.}

\begin{thm}{Lemma}\label{lem:mid} 
Let $ A$ be an Alexandrov space with curvature $\ge 0$
and $e_1, e_2,x\in A$ be two extremal points.
Assume $z$ is the midpoint of a shortest path $[e_1x]$ in $ A$.
Then 
$$|e_1z|\le|e_2z|.$$ 
Moreover, if $|e_1z|=|e_2z|$ then 
$$\angle e_1ze_2=\~\angle e_1ze_2, \quad \angle e_2zx=\~\angle e_2zx$$
and there is a unique flat % subgeodesic 
triangle 
$e_1e_2x$ in $A$ with 
a
given
 median $[e_2z]$
 (flat triangle means here the subset of $A$ isometric to the Euclidean triangle).

\end{thm}

%\begin{wrapfigure}{r}{33mm}
%\lbl[rt]{12,2;$\~e_1$}
%\lbl[lt]{135,1;$\~e_2$}
%\lbl[br]{80,34;$\~z$}
%\lbl[br]{150,58;$\~x$}
%\end{lpic}
%\end{wrapfigure}

\parit{Proof of the lemma.}
Let us assume the contrary, i.e. $|e_1z|>|e_2z|$.
Consider a comparison triangle $\~e_1\~e_2\~z$ for the chosen triangle $e_1e_2z$.
Let $\~x$ be a point on the line extension of $\~e_1\~z$
such that $|\~x\~z|=|\~z\~e_1|$.
Since $|\~e_1\~z|>|\~e_2\~z|$,
we have $\angle \~e_1\~e_2\~x>\pi/2$.
From triangle comparison, we have 
$|e_2x|\le|\~e_2\~x|$.
It follows that
$$\angle e_1e_2x\ge\~\angle e_1e_2x\ge \angle \~e_1\~e_2\~x> \pi/2.$$
In particular, $\diam\Sigma_{e_2}>\pi/2$, which is a contradiction.

In the case of equality $|e_1z|=|e_2z|$ by using the same comparison picture
as above we have $\angle \~e_1\~e_2\~x=\pi/2$.
Reasoning by contradiction, assume
$\angle e_1ze_2>\~\angle e_1ze_2$.  Then
from the
 triangle comparison we obtain $|e_2x|<|\~e_2\~x|$ and hence
$$\angle e_1e_2x\ge\~\angle e_1e_2x>\angle \~e_1\~e_2\~x= \pi/2,$$
which is a contradiction, proving angle equalities.
Now  the existence of a flat 
triangle follows from 
Lemma~\ref{lem:poly-X}.

\qeds

Now, let us introduce additional notation:
\begin{enumerate}
\item Denote by $W_i$   the set of midpoints of all geodesics $[e_ix]$ with $x\in A$.
\item Denote by $V_i$  the Voronoi domain of $e_i$;
i.e.
$$V_i=\set{x\in  A}{|e_ix|\le |e_jx|\ \text{for all}\ i}$$ 
\end{enumerate}
From Lemma~\ref{lem:mid},
we have $W_i\i V_i$ for all $i$.

Further, consider a map $\phi_i\:W_i\to  A$,
implicitly 
defined by the following relation:
 $x=\phi(z)$ if $z$ is a midpoint of a geodesic $[e_ix]$.

By triangle comparison, we have 
$$|\phi_i(z)\ \phi_i(z')|\le 2\cdot|zz'|$$
for all $z,z'\in W_i$.
In particular, the map $\phi_i$ is uniquely defined.

Hence 
$$\vol V_i\ge \vol W_i\ge \tfrac1{2^n}\cdot\vol  A.$$
Since
$$\sum_{i=1}^m\vol V_i=\vol A,$$
we get $m\le 2^n$.
\qeds

Note that from the proof we immediately get the following:

\begin{thm}{Corollary}\label{cor:V=W}
Let $\mathcal A$, $n$, $m$, $V_i$ and $W_i$ be as in the proof of Theorem~\ref{thm:2^n}.
If $m=2^n$ then $W_i=V_i$  
and
$\vol V_i= \tfrac1{2^n}\cdot\vol \mathcal A$
for all $i$.
\end{thm}

%%%%%%%%%%%%%%%%%%%%%%%%%%%%%%%%%%%%%%%%%%%%%%%%%%%%%%%%%%%%%%%%%%%%%%%%%%%

\section{Preliminary statements}

In this section we prove a number of technical statements 
needed in the proof of Theorem~\ref{main1}.
%Further we denote by $\mathcal A$ an Alexandrov space.
\subsection{Flat slices in Alexandrov space}

\begin{thm}{Lemma}\label{lem:poly-X}
Let $A$ be an $n$-dimensional Alexandrov space with nonnegative curvature and $[p x_1],[p x_2],\dots,[p x_k]$ be geodesics in $A$.

Assume that 
$$\angle x_ipx_j
=
\widetilde \angle x_i p x_j$$ 
for all $i,j$ and that all directions $\uparrow_{[p x_i]}$ lie in a subcone $E$ of $T_p A$
which is isometric to a convex cone in the
  Euclidean space.

Then all geodesics $[p x_i]$ lie in a subset of $ A$ which is isometric to a convex polyhedron in
the  Euclidean space.
\end{thm}

\parit{Proof.}
Set $\~x_i=\log_px_i\in T_p$ 
and $\~p=\log_p p$ ($\~p$ is the vertex of $T_p$).
Clearly 
\begin{itemize}
\item $\~p,\~x_1,\dots,\~x_k\in E$.
\item $|px_i|=|\~p\~x_i|$ for each $i$;
\item $|x_ix_j|=|\~x_i\~x_j|$ for all $i,j$
\end{itemize}
Since $E$ is Euclidean, 
by  Kirszbraun theorem, 
there is a short map $s: A\to E$ such that $s(p)=\~p$ and $s(x_i)=\~x_i$ for each $i$.

On the other hand the gradient exponent $\gexp_p$ is also à short map.
Thus the composition $f=s\circ\gexp_p$ is also short.
Clearly $f$ does not move $\~x_i$ and $\~p$.
It follows that $f$ does not move any point in $Q=\mathop{\rm Conv}(\~p,\~x_1,\dots,\~x_k)$. 
Therefore, $\gexp_p$ maps $Q$ isometrically into $ A$.
\qeds

\subsection{Affine functions}\label{def.aff}
In this section $A$ is an Alexandrov space of nonnegative curvature.
\begin{rdef}{Definition}
Let 
$\Omega\subset \mathcal A$ be an open subset 
and $\lambda\in\RR$.
A locally Lipschitz function $f\:\Omega\to\RR$
is called \emph{$\lambda$ quasi-affine} if 
$$(f\circ\gamma)''(t)\equiv \lambda$$
for any unitspeed geodesic $\gamma$ in $\Omega$.
We also call $0$ quasi-affine functions  \emph{affine functions}.
\end{rdef}

For an Alexandrov space
$  A$, its subset
$  \Omega\subset{ A} $
and a function
 $f\:\Omega\to\RR$
 we denote by
$ \overline{  A}$ the doubling of $  A$, by 
$ \overline{ \Omega}\subset\overline{ A} $
 the doubling of $  \Omega$
and by $\overline f:\overline \Omega\to\R$ the tautological
extension of $f$.
 
\begin{rdef}{Definition}
We say that a $\lambda$ quasi-affine
 $f\:\Omega\to\RR$
 satisfies the boundary condition if
 $\overline f:\overline \Omega\to\R$
 is $\lambda$ quasi-affine.
   \end{rdef}

For $i\in\{1,2\}$, 
assume $f_i\:\Omega\to\RR$ to be a $\lambda_i$ quasi-affine function.
Then $f_1+f_2$ is  $(\lambda_1+\lambda_2)$ quasi-affine.
Also for any real constant $c$,
 $c\cdot f_1$ is $(c\cdot\lambda_1)$ quasi-affine.

\subsubsection
{Cones and splittings}

For the proof of Proposition~\ref{prop:affine-R-factor},
 Proposition~\ref{prop:cone-structure}
and Lemma~\ref{cl:af.grad},
we refer to  \cite{a-b}.
 Functions considered in this paper are
defined on the whole Alexandrov space,
but the proof works also for  our local  
case.  It suffices to
  note that every
shortest path between points in $B_{r/4}(p)$ lies inside
$B_r(p)$.

\begin{thm}{Proposition}\label{prop:affine-R-factor}
Let  $f_1, f_2, \dots, f_k$ 
be affine functions defined on a ball
\\
 $B_r(p)\subset A$ such that the functions $1,f_1,f_2\dots, f_k$ form a linearly independent system. 
Then the ball $B_{r/4}(p)$
is isometric to an open ball in a product $\RR^k\times X$ for some metric space $X$. 
 Gradients $\nabla f_1, \nabla f_2, \dots, \nabla f_k$
 are tangent to $\R^k$ fibers. 
\end{thm}

\begin{rdef}{Definition}\label{con.negh}
A point $p\in  A$  admits a \emph{conic neighborhood} if there is an isometry from a neighborhood of $p$ to an open set in a Euclidean cone, which sends $p$ to the vertex of the cone. 
\end{rdef}

\begin{thm}{Proposition}\label{prop:cone-structure}
Suppose a ball $B_r(p)\subset A$ admits a  $1$-affine function $f$.
Then the ball $B_{r/4}$ can be isometrically
identified with  an open ball in a Euclidean cone.
Gradients $\nabla f$ are tangent to rays of the cone.
If $\nabla_p f=0$ 
we have $f=\tfrac12\dist_p^2+c$ and the ball $B_{r/4}(p)$ is a conic neighborhood of $p$.

\end{thm}

\begin{thm}{Lemma}\label{cl:af.grad}
Let $f$ be a $\lambda$ quasi-affine function
defined in some neighborhood $U\ni p$ in $A$
and $f$ satisfies the boundary condition.
Then
the tangent cone $T_pA$ splits along a line with a direction
  $\nabla_p f$ 
and $d_pf=\<\nabla_p f, \cdot\>$.

\end{thm}

\subsubsection{Dimensions of spaces of affine functions}

\begin{rdef}{Definitions}
For a  set  $F$ of affine functions
defined in some neighborhood $ U\ni p$ in $A$ we
  denote by $\#_L (F,p)$ 
the maximal number of functions in $F$, say
$f_1,\dots, f_k$, such that 
 the functions $1,f_1,f_2\dots, f_k$ form a linear independent system
 in some small ball $B_r(p)\subset U$. We note that since an affine function
 on every geodesic is determined by its initial value and its initial derivative
then $\#_L (F,p)$
 does not depend on $r$. % the choice  of    a ball $B_r(p)$.
% that are linearly independent in arbitrary small neighborhood of $p$.

For a  set  $F$ of 1-affine functions
defined in some neighborhood $U\ni p$ in $A$ we
define a set of affine functions
$F^0=\{ \sum \alpha_i f_i| f_i\in F, \alpha_i\in\R, \sum \alpha_i =0  \}$
and  define $\#_A (F,p)$ to be $\#_L (F^0,p)$.

 It follows from
Lemma~\ref{cl:af.grad}
 that the gradients of functions in $F$
 lie in a linear subspace of $T_pA$.
 Therefore we can define the following numbers:
$\#_L (\nabla F,p)$ --
the dimension of the vector subspace in $T_p\mathcal A$,
generated by the
gradients of functions in $F$ and 
$\#_A (\nabla F,p)$ -- the dimension of the affine subspace,
generated by  endpoints of these gradients.
\end{rdef}

\begin{thm}{Lemma}\label{dim.aff}
Let $F$ be a finite set of affine functions
defined 
in a  ball $B_r(p)$,
then
$\#_L (F,p)=\#_L (\nabla F,p)$.

Let $F$ be a finite set of 1-affine functions
defined 
in a  ball $B_r(p)$, then
$\#_A (F,p)=\#_A (\nabla F,p)$.

\end{thm}

\parit{Proof.}
%We choose a small ball $B_r(p)\subset U$.
It follows from
Lemma~\ref{cl:af.grad}
that the differential  of
every affine (1 quasi-affine) function is uniquely detemined by
its gradient and hence
every affine (1 quasi-affine) function $f:B_r(p)\to\R$
is determined by
 $f(p)$
and $\nabla_pf$.
%Since the set of affine  functions form linear subset  in the function space and
%the set of 1-affine functions form affine subset in the function space
Now the proof is straightforward.
\qeds
 
\begin{thm}{Corollary}\label{thm-d=d+1}
Let $F$ be a finite set of 1-affine functions
defined  in a ball $B_r(p)$.
Then the ball $B_{r/4}(p)$ can be isometrically identified with
an open ball in $\R^{\#_A (\nabla F,p)}\times  \mathcal C$,
where $\mathcal C$ is a Euclidean cone.
 If $\#_L (\nabla F,p)=\# _A(\nabla F,p) $, %=\flat ( F,x)-1$.
then the ball $B_{r/4}(p)$ is a  conic neighborhood of $p$. 
 Gradients of functions in $F$ are tangent 
 to products of  $\R^{\#_A (\nabla F,p)}$ factors
 and rays of the cone $\mathcal C$.
 
 \end{thm}

\parit{Proof.} Let's consider the set
$F^0=\{ \sum \alpha_i f_i|f_i\in F,\alpha_i\in\R, \sum \alpha_i =0  \}$.
Then $F^0$ is a set of affine functions
and
$\#_L(\nabla F_0,p)=\#_A(\nabla F,p)$,
hence by Lemma~\ref{dim.aff} and Proposition~\ref{prop:affine-R-factor}
we obtain that the ball $B_{r/4}(p)$ is
isometric to an open subset of
$\R^{\#_A (\nabla F,p)}\times X$.
Applying Proposition~\ref{prop:cone-structure}
we obtain that $X$ is isometric to an open subset in Euclidean cone.
If
$\#_L (\nabla F,p)=\# _A(\nabla F,p) $
 there are numbers $\alpha_i$, 
such that
$\sum \alpha_i=1 $
and
$\sum \alpha_i\nabla_p f_i=0 $.
Then the function
$f=\sum \alpha_i f_i $ is $1$-affine and
$\nabla_p f=0$.
Hence by Proposition~\ref{prop:cone-structure},  $f=\tfrac12\dist_p^2+c$ and 
the ball $B_{r/4}(p)$ is a conic neighborhood of $p$.
\qeds
\subsubsection{Moving Lemma}

The next lemma 
is  a technical tool for our proof of Lemma~\ref{conichelp}.
The lemma shows how
we can move a point in the domain of some collection
 of $1$-affine functions.
 Corollary~\ref{thm-d=d+1} 
makes it possible
to shift a point in a flat subset 
so that the distances  behave as Euclidean ones.

\begin{thm}{Moving Lemma}\label{lem:moving}
Let points $x, p_1,\dots, p_k\in A$ and $r>0$. 
  Suppose  that $x$ does not admit a conic neighborhood
and  the following
conditions hold: 
\begin{enumerate}[(i)]
\item The functions $f_1=\tfrac12\cdot\dist^2_{p_1}, \dots,f_k=\tfrac12\cdot\dist^2_{p_k}$ are
$1$-affine in a neighborhood $B_r(x)$.
\item $|p_1x|=|p_2x|=\dots=|p_k x|$.
\end{enumerate}

Then  there exists a unique unit vector
$v\in Span(\nabla f_1,\dots,\nabla f_k)$  such that
$\angle(v,\nabla f_1)=\dots=\angle(v,\nabla f_k)=\alpha<\pi/2$
and
 a shortest path $\gamma:[0,r/4]\to A$,
 with $\gamma(0)=x$ and  $\gamma'(0)=v$. For every
 point $y=\gamma(t)$ where
  $t\in[0,r/4]$
 we have
 the following:
\begin{enumerate}
  \item  Some small neighborhoods of $x$ and $y$ are homothetic.
 
\item $f_i(\gamma(t))=|\nabla f_i|\cos(\alpha)t+\frac{1}{2}t^2$,
in particular $|p_1y|=\dots=|p_k y|>|p_1x|$.

\item   $\angle(\gamma'(t)),\nabla f_1)=\cdots=\angle(\gamma'(t)),\nabla f_k)<\alpha$
 and
$ \# ( \nabla \{f_1,\dots,f_k\},y)=\# ( \nabla \{f_1,\dots,f_k\},x)$.
  
\item Suppose that for some $p\in A$
 the corresponding function $f_p=\tfrac12\cdot\dist_p^2$
 is $1$-affine in some neighborhood of $y$,
 $f_p(y)=f_i(y)$ and $\angle (\nabla_y f_p,\gamma'(t))\neq\angle (\nabla_y f_i,\gamma'(t))$,
 then
 $\#_A ( \{f_p,f_1,\dots,f_k\},y)=\#_A (\{f_1,\dots,f_k\},x)+1$.
 \end{enumerate}
  \end{thm}
\parit{Proof.}
 We 
apply Corollary~\ref{thm-d=d+1} and obtain
an isometric decomposition 
of $B_{r/4}(x)$ as a subset of
$\R^m\times \mathcal C $, where $m=\# _A(\nabla F,p) $
 and $\mathcal C$ is a Euclidean cone.
Vectors
$\nabla_x f_1,\dots,\nabla_xf_k$ are tangent to a 
subset of $\R^{m}\times\R_+$, namely a product of $\R^m$ and a ray in $\mathcal C$. 
We call this set a flat $(m+1)$-slice. 

 For any set $F$ of $1$-affine functions
 one of the following equalities holds: $\#_L (\nabla F,p)=\# _A(\nabla F,p)+1$, or
$\#_L (\nabla F,p)=\# _A(\nabla F,p) $.
Since $x$ does not have a conic neighborhood by Corollary~\ref{thm-d=d+1} we have
$$\#_L(\nabla \{f_1,\dots,f_k\},x)=\#_A(\nabla \{f_1,\dots,f_k\},x)+1 .$$
 Hence
 there exists a unique unit vector $v\in Span(\nabla f_1,\dots,\nabla f_k)$  such that
$\angle(v,\nabla f_1)=\dots=\angle(v,\nabla f_k)=\alpha<\pi/2$.
Then there exists   a shortest path $\gamma:[0,r/4]\to A$,
 with $\gamma(0)=x$,  $\gamma'(0)=v$ in
 our flat $(m+1)$-slice.
  Properties 1-3 follow  from the Euclidean structure.
  
   We show (4)
   arguing by contradiction. Suppose the conclusion of (4)
   does not hold, then
   $$\#_A(\nabla \{f_1,\dots,f_k,f_p\},x)=\#_A(\nabla \{f_1,\dots,f_k\},x).$$
 Hence  $\nabla_yf_p$ lies in the affine hull
  the endpoints of vectors
 $\nabla_yf_1,\dots,\nabla_yf_k$.
      We also know that 
  $|\nabla_yf_p|=|\nabla_yf_1|=|\nabla_yf_2|=\cdots=|\nabla_yf_k|$
and  
$\angle(\gamma'(t)),\nabla_y f_1)=\angle(\gamma'(t)),\nabla_y f_2)=\cdots=\angle(\gamma'(t)),\nabla_y f_k)$.  It follows that $\angle (\nabla_y f_p, \gamma'(t))=\angle (\nabla_y f_i,\gamma'(t))$, this is a contradiction.
   
    \qeds

\subsubsection{Volume evolution for a gradient flow}

Given a semiconcave function $f\: A\to\RR$,
we denote by $\Phi_f^t\: A\to A$ the corresponding gradient flow for  a time $t$.

\begin{thm}{Theorem}\label{thm:affine=concave+vol}
Let $f$ be a $\lambda$-concave function and $\Omega\subset  A$  an open set.
Then for every $t>0$, we have
$$\vol \Phi_f^t(\Omega)
\le 
\exp({n\cdot\lambda\cdot t})
\cdot
\vol \Omega.$$
Moreover if the equality holds for some $t>0$,
then $f$ is $\lambda$-affine in $\Omega$ and satisfies the boundary
condition.
\end{thm}

\parit{Proof.}
Here  $\gamma'_-$
denotes the velocity of a curve $\gamma$ if we go backwards.

$\lambda$-concavity of  $f$ means that
$$ d_pf(\gamma'(a))+d_qf(\gamma'_-(b))\ge -\lambda |pq|$$
 for every
unit speed shortest path $\gamma$ in $\Omega$
between $p$ and $q$.
To prove 
that $f$ is $\lambda$-affine
it suffices to show that this inequality turns into an equality.
We consider gradient curves $p(t)$ and $q(t)$
and let $l$ be the distance function $l(t)=|p(t)q(t)|$.

By the first variation formula
 $$l'(t)\le -(\<\gamma'(a), \nabla_p f\>+\<\gamma'_-(b), \nabla_q f\>).$$
 By definition of gradient for every  point $x$ and $v\in T_x\mathcal A$
 we have
  $ \<v, \nabla_x f\>\ge d_xf(w)$.
  Thus
  $$l'(t)\le \lambda |pq|,$$
   and applying
Proposition~\ref {prop:vol-preserv}   we obtain the  
required
volume inequality. In the case when this
inequality becomes an
equality we have that
  $l'(t)= \lambda |pq|$.
  Hence 
  $$d_pf(\gamma'(a))=\<\gamma'(a), \nabla_p f\>,\quad d_qf(\gamma'_-(b))
  =\<\gamma'_-(b), \nabla_q f\>$$
  and $\lambda$-quasi-affinity  follows.

To prove the boundary condition 
it is  enough  to check the $1$-quasi-affinity on every
shortest path 
$ \gamma:[-h,h]\to\overline{\Omega}$
intersecting  $ \overline{\partial  A}$  only once at a point
 $x=\gamma(0)\in\overline{\partial  A}$.
 Clearly,  it  suffices  to prove that
 $d_x\overline f(-\gamma'(0))=-d_x\overline f(\gamma'(0)).$
 
By above for every $x\in A\cap\Omega$
we have $d_xf=\< \nabla_xf,\cdot \>$ and
hence 
the tangent cone $T_xA$ splits along a line with a direction
  $\nabla_x f$. Then
for every $x\in\partial A\cap\Omega$  both vectors $\nabla_xf, -\nabla_xf$
lie in $\partial T_xA$
and are glued with themselves under doubling.
Hence the tangent cone of the doubling
  $T_x\overline\Omega$ also splits along
  a line with a direction
  $\nabla_x \overline f$. Thus $\angle(-\gamma'(0),\nabla f)=\pi-\angle(\gamma'(0),\nabla f)$
 and $d_x\overline f(\gamma'_-(0))=-d_x\overline f(\gamma'(0))$.

    \qeds

\subsection{Polyhedral spaces.}

\begin{thm}{Definition}
A metric on a simplicial complex $\mathcal S$ is called \emph{polyhedral} if each simplex in $\mathcal S$ is isometric to a simplex in a Euclidean space.

A metric space $\mathcal P$ is said to be \emph{polyhedral space} if it is isometric to a simplicial complex with a polyhedral metric.
\end{thm}

For the proof of
Proposition~\ref{prop:orbi-char} we need the following definition.

\begin{thm}{Definition}
A metric on a simplicial complex $\mathcal S$ is said to be \emph{spherically polyhedral} if each simplex in $\mathcal S$ is isometric to a simplex in the unit sphere in $\R^n$.

A metric space $\mathcal P$ is said to be
 \emph{spherically polyhedral space} if it is isometric to a simplicial complex with a polyhedral metric.
\end{thm}

The proof of the following characterization of polyhedral spaces 
can be found in
\cite{LP}.

\begin{thm}{Proposition}\label{prop:poly-char}
Let $X$ be a 
compact length space.
Assume that each point $x\in X$ has a conic neighborhood.
Then $X$ is a  polyhedral space.
\end{thm}

\subsection{Orbifolds.}
It is known that  for any orbifold 
that can be equipped with a metric of
constant curvature
the universal branched cover
 is a manifold.
The following proposition is 
colloquially known but we did not find appropriate reference.
This proposition 
characterizes  
$\R^n$-quotient spaces or equivalently flat orbifolds among all polyhedral
spaces.

\begin{thm}{Proposition}\label{prop:orbi-char}
A polyhedral space $P=(\mathcal S,d)$ is isometric to a quotient space $\RR^n/\Gamma$,
for a discrete action by isometries $\Gamma\acts\RR^n$
if and only if: 
\begin{enumerate}
 \item The simplicial complex $\mathcal S$ of $P$ is an $n$-dimensional pseudomanifold;
i.e.
$\mathcal S$ is connected;
any simplex in $\mathcal S$ is a face of a simplex of dimension $n$;
the link of every simplex of dimension $\le n-2$ is connected;
 every simplex of dimension $n-1$ belongs to at most two simplexes of dimension $n$.
\item For any point $x$ on a face $F$ of codimension $2$ in $P$,
the normal cone $N_xF$ of $F$ at $x$ is isometric to a quotient of $\RR^2$ by a subgroup of rotations.
Namely, $N_xF$ is  isometric to a cone over $S^1$ of length $2\cdot\pi/k$ 
or to a cone over an interval of length $\pi/k$ for some $k\in\NN$.
\end{enumerate}
\end{thm}

\parit{Proof.}
The "only if" part is obvious.
To prove the "if" part
 it is sufficient to check that $P$ is an orbifold, i.e.
 for any point $x$  in $P$
the tangent space is of the form $\RR^n/\Gamma$.
 
It is convenient to prove 
the same statement as in our Proposition
for a  spherical polyhedral space in place of polyhedral space by
and for $S^n/\Gamma$ in place of
 $\RR^n/\Gamma$.  
 So let us say that a space is 'good'  if it is polyhedral or spherical polyhedral space and
  possesses (1) and (2). 
 
 We prove by inverse induction on dimension that every
 'good' space is isometric to $\RR^n/\Gamma$ or $S^n/\Gamma$.
%(note again that it is sufficient to prove that 'good' space is an orbifold).
 The base $k=2$ follows because of Condition 2.
 Suppose any 'good' space of dimension $k-1$
 is isometric to $\RR^{k-1}/\Gamma$ or $S^{k-1}/\Gamma$.
 Then for any $k$-dimensional 'good' space $P$ and any point 
 $x\in P$ the
 unit tangent space $\Sigma_xP$ is a spherical polyhedral space
 which inherits properties (1) and (2) and hence is 'good'.
 Hence by the induction   hypothesis
  $\Sigma_xP=S^k/\Gamma$ and $P$ is an orbifold.
  This proves the induction step.
 
  \qeds

\subsection{Volume preserving + 1-Lipschitz  = isometry}
The proof of the following fact can be found in
\cite{nanli}.

\begin{thm}{Proposition}\label{prop:vol-preserv}
Let $\mathcal X$ and $\mathcal Y$ be $m$-dimensional Alexandrov spaces,
$\Omega\subset \mathcal X\backslash \partial \mathcal X$  an open set
and $f\:\Omega\to Y$  a 1-Lipschitz volume preserving map.
Then $f$ is a locally distance preserving;
i.e., for every point $x\in \Omega$ there is a neighborhood $\Omega_x\ni x$
such that the restriction $f|\Omega_x$ is a distance preserving map.
\end{thm}

\section{Any $n$-box is a polyhedral space.}

\begin{rdef}{Notation}
In what follows we denote by $\mathcal A$  an $n$-box.

We keep the notation for $e_i$, $V_i$, $W_i$ and $\phi_i$ for all $i\in\{1,2,\dots,2^n\}$ from Subsection~\ref{subsec:2^n}.
According to  Corollary~\ref{cor:V=W}, $V_i=W_i$ for all $i$.

Denote by $\mathfrak C_i$ the cutlocus of $e_i$;
i.e. the set of points $z\in \mathcal A\backslash\{e_i\}$ which do not lie in the interior of 
every
% geodesic
 shortest path $[e_ix]$. 
\end{rdef}

In this section we prove the following  result:

\begin{thm}{Theorem}\label{thm:polih}
Every $n$-box is a polyhedral space.
\end{thm}

 Proposition~\ref{prop:poly-char} implies that
it is sufficient to prove the following lemma:

\begin{thm}{Key Lemma}\label{lem:cone}
Every point $x\in \mathcal A$ has a conic neighborhood.
\end{thm}

\begin{thm}{Proposition}\label{prop:lochom}
Each function $f_i=\tfrac12\cdot\dist^2_{e_i}$ is $1$-quasi-affine
and satisfies the boundary condition
 in $\mathcal A\setminus \mathfrak C_i$.
\end{thm}

\parit{Proof.}
It is sufficient to note that 
the restriction $\Phi^{\ln 2}_{f_i}|_{W_i}$
coincides with
$\phi_i$. 
Then  from 
%Proposition~\ref{prop:vol-preserv} 
Corollary~\ref{cor:V=W}
and Theorem~\ref{thm:affine=concave+vol} it
follows that $f_i$ is $1$-quasi-affine and satisfies the
boundary condition.

 \qeds

%Applying propositions \ref{prop:lochom} and \ref{prop:cone-nbhd}, 
%we get
 
%\begin{thm}{Corollary}\label{cor:loc.cone} 
%Let $p\in \mathcal A\setminus \mathfrak C_i$.
%Then the necessury unique geodesic $[e_i p]$ has a neigborhood which admits an isometry to an open set in a cone. 

%Moreover one can choose an isometry which sends $e_i$ to the vertex of the cone.
%\end{thm}

\parit{Proof of the Key Lemma.}
It follows 
from Proposition~\ref{prop:lochom}, 
that for every point $x\in V_i$ the function $f_i$ 
is $1$-affine in a neighborhood of $x$.  Therefore,
we can define
for 
a given  point $x\in\mathcal A$
 an
 index set $J_x\subset\{1,\dots,n\}$ and a positive integer
 $\#(x)$ as follows:
  $$J_x= \{i\in\{1,\dots,n\}| x\in V_i\};$$
$$\#(x)=\#_A\{f_i|i\in J_x\}.$$

%to be the maximal number of elements in $Q\subset J_x$  such that the functions $\set{f_i}{i\in Q}$ are linearly independed in arbitrary small neighborhood of $x$.

According to Lemma~\ref{dim.aff} and Corollary~\ref{thm-d=d+1} we have that
 $\#(x)\le n$ for every $x\in \mathcal A$.
Moreover if $\#(x)= n$ then $x$ has a flat neighborhood.

The main technical point of the proof of the Key Lemma is the following:
\begin{thm}{Lemma}\label{conichelp}
Assume a point $x\in \mathcal A$ has no conic neighborhood.
Then there is a point $x'\in \mathcal A$ such that
a neighborhood of $x'$ is homothetic to a neighborhood of $x$ and 
$\#(x')>\#(x)$.
\end{thm}

We prove this lemma in Subsection~\ref{profconhelp}.

Now to
prove the Key Lemma ~\ref{lem:cone} we argue by contradiction. 
Let us assume the contrary,
i.e., there is a point $x\in \mathcal A$ which has no conic neighborhood.

Applying Lemma~\ref{conichelp} for $x_0=x$, we get a point $x_1$ with a neighborhood homothetic to a neighborhood of $x_0$ and $\#(x_1)\ge\#(x_0)+1$.
In particular, $x_1$ does not admit a conic neighborhood.

Therefore we can apply Lemma~\ref{conichelp} $(n+1)$ times to get a
point $ x_{n+1}\in\mathcal A$
such that $\#(x_{n+1})\ge n+1$.
We arrive to a contradiction since $\#(z)\le n$ for any $z\in \mathcal A$.
\qeds

\subsection{Proof of Lemma~\ref{conichelp}}  \label{profconhelp}

For each $i$ the sets
$V_i$ and $\mathfrak C_i$ are closed and  
disjoint.
Hence Proposition~\ref{prop:lochom} implies that
 there exists
$r_0>0$ such that for every $i$ and $x\in V_i$ the  function $\frac{1}{2}\dist_{e_i}^2 $ 
is $1$-quasi-affine in $B_{4r_0}(x)$.

%Let note that in our notations $\#(x)=\flat(\{f_i|i\in J_x\},x)$.

Now we fix $x\in \mathcal A$ and suppose that $x$ does not have a conic neighborhood.
We apply Lemma~\ref{lem:moving} for $x$ and 
 $ \{f_i|i\in J_x\}$. We can shift $x$ 
 equidistantly from points $e_i$ for $ i\in J_x$
so that
the points still lie in all $V_i$ for $i\in J_x$.
We continue until  we meet a domain $V_j$ for some $j \not\in J_x$.
Let us formulate the exact statement:
 
Let  $\gamma_0:[0,r_0]\to A$ be the shortest path from Lemma~\ref{lem:moving}.
We have a dichotomy:
 
\begin{enumerate}
\item\label{one} There exists a minimal value $t_0\in(0, r_0]$, such that
$\gamma_0(t_0)\in V_j$ for some $j_0\not\in J_x$.
Set $y=\gamma_0(t_0)\in V_{j_0}$ and $f_i=\frac{1}{2}\dist_{e_{i}}^2 $,
for $i\in J_x\cup \{j_0\}$.
We have the angle 
   inequality
 $\angle (\nabla_y f_{j_0},\gamma'(t))>\angle (\nabla_y f_i,\gamma'(t))$
(indeed, otherwise we would have that
 $f_{j_0}(t_0-\epsilon)\le f_i(t_0-\epsilon)$
 for sufficiently small $\epsilon>0$, this would contradict the choice of 
 $t_0$).  
Thus we can apply Moving Lemma~\ref{lem:moving}  (4)
 with $p:=e_{j_0}$, $f_p=f_{j_0}=\frac{1}{2}\dist_{e_{j_0}}^2 $. Then  some small neighborhoods of $x$ and $y$ are homothetic
 and
\begin{align*}
 \#(y)
&\ge
\# _A( \{f_i|i\in J_x\}\cup\{f_{j_0}\},y)=
\\
&=\# _A(\{f_i|i\in J_x\},x)+1=
\\
&=\#(x)+1.
\end{align*}

\item The shortest path $\gamma_0([0,r_0])$ does not intersect any $V_j$  for
  $j\notin J_x$.
 
In this case we apply Moving  Lemma (\ref{lem:moving})
recursively for $x_1=\gamma_0(r_0)$  and so on.
After $k$ iteration we have an estimate 
$f_i(x_k)>(|\nabla_x f_i|\cos(\alpha_0)r_0)\cdot k$, $i\in J_x$
    where $\alpha_0=\angle(\nabla_x f_i,\gamma_0'(0))$.  
The diameter of $A$ is finite, 
    therefore after
    finitely many steps we arrive to  Case~\ref{one}.\qeds
\end{enumerate}

\section{$n$-boxes are flat orbifolds}

In this section we finish the proof of Theorem~\ref{main1}.

Note that according to Theorem~\ref{thm:polih} and Proposition~\ref{prop:orbi-char},
it  suffices to show the following:

\begin{thm}{Theorem}\label{cone:::pi}
Let an $n$-dimensional polyhedral space $\mathcal A$ be a box.
Then the normal cone for each face of codimension $2$ in $\mathcal A$ is isometric to one of the following spaces:
$\RR^2$, $\RR_+\times\RR$, $\RR_+\times\RR_+$ or a cone over a circle of  length $\pi$.
\end{thm}
The proof of this theorem is in Subsection~\ref{provecone}.

Let $\mathcal A$ be an $n$-box.
We keep the same notation as above: 
$e_i$ denote extremal points of $\mathcal A$,
$V_i$ the corresponding Voronoi domain,
$\mathfrak C_i$ the cut locus of $e_i$;
$i\in\{1,2,\dots,2^n\}$. 
A minimizing geodesic $[e_i e_j]$ between two extremal points is called an \emph{edge}.

Let $p\in \mathcal A$ be a point
which lies on a face of codimension $2$;
i.e., $T_p\mathcal A= \RR^{m-2}\times L$,
where $L$ denotes a 2-dimensional cone  containing no lines.
Take the set of all points in $\mathcal A$ with tangent cone isometric to $T_p\mathcal A$;
we call
its closure $H$  \emph{hyperedge} (we name it this way since $H$ has codimension 2 in $\mathcal A$).

%??? We say that edge $[e_ie_j]$ is transversal to hyperedge $H$ if there is a geodesic $[e_ix]$ arbitrary close to $[e_ie_j]$ such that $x\in H$
%and $\l]e_ix\r[$ does not contain points in $H$. 

Here are simplest properties of hyperedges of $n$-boxes:

\begin{thm}{Lemma}

(1) Any hyperedge contains at least one   vertex $e_i$.
 
(2) If a vertex $e_j\notin H$  
then $H\subset\mathfrak C_j$.

\end{thm}
\parit{Proof.}

(1) Indeed, take a point $x\in relint (H)$
then for some $i\in\{1,\dots, 2^n\}$
$x\in V_i$. Then $e_i\in H$.

 (2) For any point $x\in relint (H)$
we have $x\in\mathfrak C_j$.
Hence $H\subset\mathfrak C_j$.   
  \qeds
\subsection{Proof of Theorem~\ref{cone:::pi}}\label{provecone}

\begin{rdef}{Definition} 
Let $\mathcal A$ be a box
and $H\subset \mathcal A$   a hyperedge.  
We say that a vertex $e_i\in \mathcal A$
\emph{pushes 
 $H$  in a vertex   $e_j\in H$}
 if there exists a flat $(n-2)$-dimensional simplex $\Delta\subset H$,
 such that $e_j\in\Delta\subset \mathfrak C_i$.
 \end{rdef}
\begin{rdef}{Definition} \label{separates}
Let $\mathcal A$ be a box
and $H\subset \mathcal A$   a hyperedge. We say that $H$  
\emph{separates a vertex $e_i\in \mathcal A$ from a vertex $e_j\in H$}
 if there exists a flat $(n-2)$-dimensional simplex $\Delta\subset H$,
 such that $e_j\in\Delta\subset \mathfrak C_i$ and
$$relint(\varphi_i^{-1}(\triangle))\cap V_k=\emptyset\ \ \text{for every}\ \  k\neq i,j.$$
 \end{rdef}

To prove Theorem~\ref{cone:::pi},
we need the following lemma:
\begin{thm}{Lemma}\label{lem:existsubdivide}
Let $\mathcal A$ be a box 
and $H$ a hyperedge.
Then there are vertices $e_i\in\mathcal A$ and $e_j\in H$
such that
 $H$
 separates $e_i$ from $e_j$.
  \end{thm} 
The proof of this lemma is in
Subsection~\ref{proofinduction}.
Now let us show how Theorem~\ref{cone:::pi} follows from Lemma~\ref{lem:existsubdivide}.

\parit{Proof of Theorem~\ref{cone:::pi}.}
Let us introduce some notation:
\begin{itemize}
\item $K_i$ denotes the completion of $\mathcal A\setminus \mathfrak C_i$ equipped with 
an intrinsic metric.

\item Clearly $K_i$ is isometric to $2{\cdot}V_i$. 
Denote by $\psi_i:g_i^{-1}(V_i)\to K_i$ the homothety  centered at $e_i$ and
with  coefficient $2$.

\item $g_i:K_i\to \mathcal A$ is the corresponding gluing map (which is piecewise linear).
\end{itemize}
Note that in these notation we have $ g_i\circ\psi_i\circ g_i^{-1}=\phi_i $.

To prove Theorem~\ref{cone:::pi}
we take a hyperedge $H$ containing a given $(n-2)$-dimensional face,
 apply Lemma~\ref{lem:existsubdivide}
and  obtain  that $H$ separates some vertices $e_i\in\mathcal A$, $e_j\in H$. Let $\Delta$
be from Definition~\ref{separates}.
It is sufficient now
to prove that for some  point $x\in relint(\Delta)$
the normal cone to $H$ at this point
is one of the 4 cones described in the theorem. 
We can assume that $\Delta$ is sufficiently small
so that $g_i^{-1}(relint(\Delta))$ are disjoint isometric copies of $relint(\Delta)$.
By $\Delta_1,\dots,\Delta_l$
we denote closures of its preimages
and by $e_j^1\in\Delta_1,\dots,e_j^l\in\Delta_l$
the corresponding preimages of $e_j$.

The next lemma describes a possible structure of
the tangent space
of a point in the preimage $g_i^{-1}(int(\Delta))$.

\begin{thm}{Lemma}\label{tangR}
Using our notation, 
let a point $x\in g_i^{-1}(int(\Delta))\subset\partial K_i$.
Then there are 2 possibilities: 

(1) if $\psi_i^{-1}(x)\notin \partial K_i$ then $T_xK_i=\R^{n-1}\times \R_+ $;

(2) if $\psi_i^{-1}(x)\in \partial K_i$ then
$T_xK_i=\R^{n-2}\times \R_+\times \R_+ $.
\end{thm}

\parit{Proof.}
We can assume that $x\in\Delta_1$.
For  $y=\psi_i^{-1}(x)$ we know, that $T_yK_i$
contains an isometric copy of $\R^{n-2}\times\R$.
Hence 
$T_yK_i=\R^n$ or $T_yK_i=\R^{n-1}\times\R_+$.
We know also that $\psi_i^{-1}(\Delta_1)$ is a flat $(n-2)$
simplex equidistant 
from $e_i$ and $e_j^1$
with midpoint $\psi_i^{-1}(e_j^1)$ as a  vertex.
In a small neighborhood $U$ of $y$
we  have that
$$g_i^{-1}(V_i)\cap U
=\set{z\in U}{|ze_i|\le|ze_j^1|}.$$

It follows that $T_y(g_i^{-1}(V_i))$ can be presented
as one  part of
perpendicular bisection of 
$T_yK_i$ w.r.t. $e_ie_j^1$.
Thus we have:
\begin{itemize}
\item $T_y(g_i^{-1}(V_i))=\R^{n-1}\times \R_+ $
 if
$T_yK_i=\R^n$;
\item $T_y(g_i^{-1}(V_i))=\R^{n-2}\times \R_+\times \R_+ $ if
$T_yK_i=\R^{n-1}\times\R_+$.\end{itemize}
It remains to note that $T_y(g_i^{-1}(V_i))$
is isometric to 
$T_xK_i$. \qeds

\begin{thm}{Lemma}\label{boundary}
 For a 
point $x\in g^{-1}(\mathfrak C_i)\subset\partial K_i$ the condition
$\psi_i^{-1}(x)\in \partial K_i$ implies $g_i(x)\in\partial \mathcal A$.
\end{thm}

\parit{Proof.}
This follows from the fact that our space is polyhedral and
$g_i(\partial K_i)\setminus\mathfrak C_i\subset \partial\mathcal A$.
\qeds

We can consider the space $K_i$ as the result of a cutting off the polyhedral space $\mathcal A$
along $(n-1)$-polyhedral subspace $\mathfrak{C}_i$.
The map $g_i$ glues $\mathcal A$ back from $K_i$.
Then if the point $x\in \mathfrak{C}_i$
has $l$ preimages  $x_1,\dots,x_l\in K_i$ under $g_i$,
its tangent space $T_x$ can be glued out from the tangent spaces $T_{x_1},\dots,T_{x_l}$.
We write this:
$$T_x=T_{x_1}\sqcup\dots\sqcup T_{x_l},$$
the gluing maps are $d_{x_1}g_i:T_{x_1}\to T_x,\dots,d_{x_l}g_i:T_{x_l}\to T_x$.

Fix $x$ and let $g_i^{-1}(x)=\{x_1,\dots,x_l\}\subset K_i$.
Then there are two possibilities:

\begin{enumerate}
\item $\Delta\subset\partial \mathcal A$

\begin{enumerate}
\item for some $1\le k_0\le l$ the point $\psi^{-1}(x_{k_0})\notin\partial K_i$.
Then by Lemma~\ref{tangR} $T_{x_k}=\R^{n-1}\times \R_+$, $l=1$ and
$T_{x}=\R^{n-1}\times \R_+$.

\item for all $k\in\{1,\dots,l\}$ points $\psi^{-1}(x_k)\in\partial K_i$.
Then $T_{x_k}=\R^{n-2}\times \R_+\times \R_+$
This is only possible
if $l=1$ and   $T_{x_k}=\R^{n-2}\times \R_+\times \R_+$ or
$l=2$ and
$T_{x}\mathcal A=\R^{n-1}\times \R_+$.
\end{enumerate}

\item 
$int(\Delta)\cap\partial \mathcal A=\emptyset$, in this case Lemma~\ref{boundary}
implies that
 for all $k\in\{1,\dots,l\}$,  $\psi^{-1}(x_k)\notin\partial K_i$
and by Lemma~\ref{tangR}
  $T_{x_k}=\R^{n-1}\times \R_+$. This only possible
if  $l=1$ and $T_x\mathcal A=\R^{m-2}\times L$, where
$L$ is a cone over $S^1$ of length $\pi$
 or $l=2$ and
$T_x\mathcal A=\R^n$.

\end{enumerate}
This completes the proof of
Theorem~\ref{cone:::pi}.
\qeds

\subsection{Proof of Lemma~\ref{lem:existsubdivide}.}\label{proofinduction}
Let us note  that if there is a vertex, say $e_1\notin H$,
then the proof would be much simpler. 
It would be sufficient to take the shortest edge
between vertices   in $H$ and outside $H$. So the difficulty is if there is no such a vertex.

To find vertices separated by $H$
we start with  Lemma~\ref{transedge} 
to find a pair of vertices $e_i$, $e_j$
such that  $e_i$ presses down $H$ at  $e_j$.
Then we can decrease the distance $|e_ie_j|$ between points
with  the same property using
Lemma~\ref{cl:shorteredge} until
we find a pair of
vertices, such that $H$ separates one from the other.

\begin{thm}{Lemma}\label{transedge}
Let $\mathcal A$ be a box.  Then for any hyperedge $H\subset\mathcal A$ 
 there are vertices $e_i\in\mathcal A$ and $e_j\in H$
such that $e_i$ pushes
 $H$ in
  $e_j$.

\end{thm}

\parit{Proof.}
Suppose there exists at least one
vertex $e_i\notin H$, then  $e_i$ pushes
$H$ in every vertex $e_j\in H$.
Otherwise consider any flat $n$-simplex
with vertexes in $\{e_1,\dots, e_{2^n}\}$ say $\triangle_{e_{i_0},\dots,e_{i_n}}$.
The existence of such a simplex
can be proved by using the same construction
as in the proof of \ref{profconhelp}:
moving out from vertexes we can find a point $x\in\mathcal A$ with
$\#(x)=n$ and from \ref{lem:poly-X} it follows that corresponding
$n+1$ vertexes form flat $n$ simplex.
Since the codimension of $H$ is $2$,
 one of the vertexes
$e_{i_1}$, $e_{i_2},\dots,e_{i_n}$
has to push $H$ in $e_{i_0}$.
\qeds

\begin{thm}{Lemma}\label{cl:shorteredge}
Let $e_i$ and $e_j$ be two vertices and
$H$  a hyperedge in an $n$-box $\mathcal A$
and $e_j\in H$.
Assume $e_i$ pushes  $H$ in $e_j$ 
but $H$ does not separate
$e_i$ from $e_j$. 
Then there is $k\neq i,j$ such that
$$\max\{|e_ke_i|,|e_ke_j|\}<|e_ie_j|$$ 
and one of the following holds:

\begin{itemize}
\item $e_k$ pushes $H$ in $e_j$;
\item $e_k\in H$ and $e_i$ pushes $H$ in $e_k$.
\end{itemize}

\end{thm}

To prove Lemma~\ref{cl:shorteredge}
we  need the following:

\begin{thm}{Subemma}\label{lem:sh.pth}
For any vertices $e_i$, $e_k$ and  a point $x\in V_i\cap V_k$
there is a shortest path $[\phi_i(x)e_k]$ inside $\mathfrak C_i$.

\end{thm}

\parit{Proof.}
By Lemma~\ref{lem:mid}
 there is a flat triangle $e_ie_k\phi_i(x)$ with median 
$[xe_k]$ and right angle at $e_k$.
If some point of the edge
$[\phi_i(x)e_k]$
of this triangle does not lie in
$\mathfrak C_i$ then we would have
$\diam\Sigma_{e_k}>\pi/2$, contradiction.
\qeds

\parit{Proof of Lemma~\ref{cl:shorteredge}.}
In conditions of our lemma there exists an
$(n-2)$-simplex $\Delta$ with a vertex $m\in\phi_i^{-1}(e_j)$
such that $\phi_i(\Delta)\subset H$
and
 $\Delta\subset V_i\cap V_k$ for some $k\neq i,j$. 
Then by Lemma~\ref{lem:mid}
 there is a flat triangle $e_ie_je_k$ with median 
$[e_km]$ and right angle in $e_k$.
Then 
$$\max\{|e_ke_i|,|e_ke_j|\}<|e_ie_j|.$$

Now if $e_k$ prushes $H$ in $e_j$ the proof is completed.
Suppose contrary. We can assume that $relint(\phi_i(\Delta))\subset \mathcal A\setminus\mathfrak C_k$.
By Sublemma~\ref{lem:sh.pth}
for every point $y\in\phi_i(\Delta)$
there is a shortest path $[ye_k]$ inside $\mathfrak C_i$, if in addition
$y\notin \mathfrak C_k$  then $[ye_k]\subset H$.
Then points of all such shortest paths for $y\in relint(\phi_i(\Delta))$ form 
an $(n-2)$-dimensional subset of $H$.
In particular $e_k\in H$ and $e_j$ presses down $H$ at $e_k$.
\qeds

\section{The structure of the action of the orbifold group of $n$-box.}
Now we are in position to
 prove 
Theorem~\ref{main2}.

In what follows we
assume that $\Gamma\acts\RR^n$ is a discrete
cocompact
 action by isometries and 
$\Pi: \RR^n\to\RR^n/\Gamma$ denotes the projection.
Let us denote by $\mathcal E$ the set of % isolated fixed points of some subgroups of $\Gamma$.
singular points for $\Gamma\acts\RR^n$.

It follows from
Proposition~\ref{prop:affine+extreme}
that $\RR^n/\Gamma$ is a box iff the number of $\Gamma$-orbits in $\mathcal E$ is $2^n$.
This implies in particular the first part of Theorem~\ref{main2}. 
We reduce the second part of Theorem~\ref{main2} to three propositions below in this section.
 To formulate the propositions we need some definitions and notation.

%In what follows  $\mathcal A$ denotes a box and $\Gamma\acts\RR^n$ a correspondent group action.
\begin{rdef}{Notation}

For any $x\in \mathcal E$ we  denote by $V_x$ its Voronoy cell w.r.t. $\mathcal E$;
i.e.
$$V_x=\set{z\in\RR^n}{|z-x|\le|z-y|\ \text{for every}\ y\in \mathcal E}.$$

Given $x\in\RR^n$, we denote by $\Gamma^\#_x\subset O(n)$ the
action 
 of the 
stabilizer $\Gamma_x$ on the vector space $\R^n$.
\end{rdef}

\begin{rdef}{Definitions}
 We say that  an action $\Gamma\acts\RR^n$ has a \emph{reflection property}
if  $\mathcal E\neq\emptyset$ and
for any adjacent $x,y\in \mathcal E$ (i.e.  $\dim(V_x\cap V_y)=n-1$)
the stabilizer $\Gamma_x$ can only fix  or reflect the point $y$:
$\Gamma_x^\#(\{\vr{xy}\})=\{\vr{xy},-\vr{xy}\}$.

We say that a discrete subset $E\subset\R^n$ is a \emph{lattice}
if there is a finite set of
\emph{generating}
 vectors $\vr{a_1},\dots,\vr{a_l}$ such that for any point $x\in E$
 we have:
$$E=\{x+k_1\cdot\vr{a_1}+
k_2\cdot\vr{a_2}+
\dots+
k_l\cdot\vr{a_l}|k_1, k_2\dots,k_l\in\ZZ \}.$$
If the dimension of the affine hull of $E$ equals $k$ we say that
$E$ is a \emph{$k$-lattice}.

We say that a group $\Gamma\acts\R^n$ \emph{reflects  generating vectors}
$\vr{a_1},\dots,\vr{a_l}$
if for any $x\in E$ and $i=1,\dots,l$ we have that
$\Gamma_x^\#(\{\vr{a_i}\})=\{\vr{a_i},-\vr{a_i}\}$.
\end{rdef}

The second part of Theorem~\ref{main2}
follows from Theorem~\ref{main1}
and the next three Propositions.

\begin{thm}{Proposition}\label{ourLattice}
Assume that
 the number of $\Gamma$-orbits in $\mathcal E$ be $2^n$.
Then  $\Gamma\acts\RR^n$  has a reflection
  property.

\end{thm}
\parit{Proof.}
%Follows from Lemma~\ref{cl:s.edgeVor}  and Lemma~\ref{SEdgeProperties}.
See proof  in Section~\ref{sim.edge}.

\begin{thm}{Proposition}\label{latticemain}
Let an action  $\Gamma\acts\R^n$ have a reflection
 property.
Then $\mathcal E$ is an $n$-lattice.
  Moreover,
there exist $n$ generating vectors for $\mathcal E$
and $\Gamma$ reflects these generating vectors.

\end{thm}

\parit{Proof.} See proof is in Section~\ref{secLattice}.

\begin{thm}{Proposition} 
Let the number of $\Gamma$-orbits in $\mathcal E$ be $2^n$.
Suppose that 
 $\mathcal E$ is a lattice
  and
there exist $n$ generating linearly independent vectors $\vr{a_1},\dots,\vr{a_n}$
 for $\mathcal E$
such that $\Gamma$ reflects this generating vectors.

Then the action
  $\Gamma\acts\R^n$
  is affine conjugate to
  a subaction of the Coxeter group assotiated to the unit cube.% $\Delta^n\acts\RR^n$.  
  \end{thm}
\parit{Proof.}
Let us denote by $\Gamma_*$ 
 a subgroup of $\Gamma$  generated by stabilizers 
of all singular points. 
Let us denote by $2\cdot\mathcal E$ the set of vectors
$\{\sum_{i=1}^n 2\alpha_i\vr{a_i}|k_1,\dots,k_n\in \ZZ\}$.
Since $\Gamma_*$ reflects the generating set
we have that
for any $y\in\mathcal E$,
$\Gamma_*(y)\supset y+2\cdot\mathcal E$.
The set $\mathcal E$ is invariant under $\Gamma$
and the number of orbits equals $2^n$. Hence we have that
$\Gamma_*(y)=y+2\cdot\mathcal E$ for any $y\in\mathcal E$.
By the same arguments we obtain that
$\Gamma(y)=y+2\cdot\mathcal E$ for any $y\in\mathcal E$.
Therefore $\Gamma(y)=\Gamma_*(y)$, hence 
$\Gamma=\Gamma_*$.

We fix coordinates in $\R^n$: a point
$O\in\R^n$  and an orthonormal basis $e_1,\dots,e_n$.
Let
$\bigoplus^n\acts\RR^n$
be the corresponding action of the Coxeter group  of the
unit cube.
We define an affine map on the basis:
$F(x_0)=O$, for some point $x_0\in\mathcal E$ and
$F(a_i)=e_i$. 
We define
 an action  $G\acts\R^n$ by
$G=F\circ\Gamma\circ F^{-1}$, this action is 
affine conjugate to
the action $\Gamma\acts\R^n$.
We have that the integer lattice $\ZZ^n$
is the set of singular points of the action
$G\acts\R^n$,
the group $G$ is generated by stabilizers of points of $\ZZ^n$
and $G$ reflects 
 the generating set
 $e_1,\dots,e_n$ of the lattice
  $\ZZ^n$.
  It follows that $G\le \bigoplus^n$. 

\section{Properties of a group action for an $n$-box.}\label{sim.edge}
In this section we prove 
Proposition~\ref{ourLattice}. 
The quotient space $\R^n/\Gamma$ is an $n$-box, we denote it by $\mathcal A$
and keep all notation for $n$-boxes,  we used before.   

We precede the proof by three lemmas.
The first two  are technical facts about Voronoy domains,
and Lemma~\ref{lem:factor.sphere} is the main
 geometric observation for our proof of Proposition~\ref{ourLattice}:

\begin{thm}{Lemma}\label{lem:liftVoronoy}
Let $M=\R^k$ or $M=S^k$
and $G\acts M$ be a discrete cocompact action by isometries.
Let us denote the quotient space $M/G$ by $M'$. Let $p$ be the projection
$M\to M'$. We fix some finite collection of points 
$s_1,\dots,s_l\in M'$ and consider Voronoy decompositions
of $M'$ and of $M$ 
w.r.t. the sets $\{s_1,\dots,s_l\}$ 
and
$p^{-1}(\{s_1,\dots,s_l\})$ correspondently.
Then
for every $i\in\{1,\dots,l\}$ and every point $s\in p^{-1}(s_i)\subset M$
the corresponding Voronoy domain $V_{s}\subset M$ can be characterized
by the following property:
 a point $y\in V_{s}$ iff
 $(y\in p^{-1}(V_{s_i}))\& (|sy|=|s_ip(y)|)$.

% we have $p(V_{e_i^0})=V_i$.
\end{thm}

\parit{Proof.} The proof is straightforward  and uses
just two properties of the projection map:  the map $p$ doesn't increase 
distances and for any $x,y\in M'$ there exist $x^0\in p^{-1}(x), y^0\in p^{-1}(y)$
such that $|x^0y^0|=|xy|$. \qeds
%any piecewise linear path in $M'$ can be lifted to the path of the same length in $M$. 

\begin{thm}{Lemma}\label{lem:Voronoyintersection}
For any two points $e_i^0\in \Pi^{-1}(e_i)$ and 
 $e_j^0\in \Pi^{-1}(e_j)$
 such that
 $V_{e_i^0}\cap V_{e_j^0}\neq \emptyset $
 the projection $\Pi$ is a distance preserving map, that is
 $|\Pi(e_i^0)\Pi(e_j^0)|=|e_i^0e_j^0|=|e_ie_j|$. 
  \end{thm}
\parit{Proof.}
Let a point $x$ lie in $V_{e_i^0}\cap V_{e_j^0}$.
By Lemma~\ref{lem:liftVoronoy} we have that 
$\Pi(x)\in V_{i}\cap V_{j}$. Then by
 Lemma~\ref{lem:mid} there exists a unique flat 
 totally geodesic triangle $e_ie_j\Pi(x)$ and the triangle
 $e_i^0e_j^0x$ is its  isometric lifting. \qeds

\begin{thm}{Lemma}\label{lem:factor.sphere}
Let $S^k$ be a $k$-dimensional sphere, $G$  a discrete subgroup
of isometries of $S^k$, $B=S^k/G$ with a projection $p:S^k\to B$ and
$diam  B\le\pi/2$.
Suppose that for a point
$v\in S^k$ the following holds:
there exists a $(k-1)$-dimensional subset $F\subset B$
such that $|p(v)x|=\pi/2$ for all $x\in F$ (further we refer to this as "$\pi/2$-property") .
Then 
 the orbit of $v$ contains exactly  two points:
  $G(v)=\{v,v^-\}$, where $v^-\in S^k$ is the diametrical point for $v$. 
\end{thm}
\parit{Proof.}
For a point $w\in S^k$,
we denote by $S_w^\bot$ the equator $\{y\in S^k||yw|=\pi/2\}$.
We consider the Voronoy decomposition of $S^k$ 
w.r.t. the set $p^{-1}(p(v))$.
Lemma~\ref{lem:liftVoronoy}
implies that
$p^{-1}(F)\cap V_{w}\subset S_w^\bot$ for any $w\in p^{-1}(p(v))$.
Since  for any $w\in p^{-1}(p(v))$
$\diam(V_w)\le\pi/2 $ and
 $\dim(p^{-1}(F))=k-1$  there are exactly two
 Voronoy domains, which are semi-spheres and the set $p^{-1}(p(v))=G(v)$
consists of two diametric points.
\qeds

\parit{Proof of Proposition~\ref{ourLattice}}

Let us  fix some adjacent points $e_i^0\in \Pi^{-1}(e_i)$,
 $e_j^0\in \Pi^{-1}(e_j)$.  Let $m^0$ be a midpoint between them.
 Let $m=\Pi(m^0)$.
  It follows from 
 Lemma~\ref{lem:Voronoyintersection}  
 that $|e_i^0e_j^0|=|e_ie_j|$ and 
 $|e_im|=|e_jm|$, hence by Corollary~\ref{cor:V=W} we 
 have that the midpoint
 $m$ lies in $V_i\cap V_j$. Then 
  applying 
   Lemma~\ref{lem:liftVoronoy}
   we obtain that $m^0\in   V_{e_i^0}\cap V_{e_j^0}$.
   It follows from the definition of adjacent vertices 
   and convexity 
   of Voronoy domains that
    there exists a flat $(n-1)$-triangle $\Delta^0$ with a vertex $m^0$ such that
     $\Delta^0\subset V_{e_i^0}\cap V_{e_j^0}$.  
 Hence there exists a flat $(n-1)$-triangle $\Delta$ with a vertex $m=\Pi(m^0)$ such that
     $\Delta\subset V_{e_i}\cap V_{e_j}$ 
     (indeed we can take a sufficiently small triangle  in $\Pi(\Delta^0)$).
      
  We consider the action of the stabilizer on the unit sphere
 $\Gamma_{e_j^0}\acts \Sigma_{e_j^0}\R^n$ and the 
 corresponding
 quotient space $\Sigma_{e_j^0}\R^n/ \Gamma_{e_j^0}=\Sigma_{e_j}\mathcal A$,
 let denote the corresponding projection by $p:\Sigma_{e_j^0}\R^n\to\Sigma_{e_j}\mathcal A$.
  To prove the reflection
 property
 we apply
 Lemma~\ref{lem:factor.sphere} 
   to this action, vector $v=\vr{e_j^0e_i^0}/|e_j^0e_i^0|$ (as a point in $\Sigma_{e_j^0}\R^n$)
  and the set $F=d_{m}\phi_i(\Sigma_{m}\Delta)$. It is not difficult to see
that the map $d_{m}\phi:T_mV_i\to T_{e_j}\mathcal A$ does not decrease
dimensions (indeed, this map is a gluing map), hence $\dim F = n-2$.  
 To verify conditions
 of Lemma~\ref{lem:factor.sphere} it remains to prove 
 $\pi/2$-property for the vector $v$ and the set $F$.
 For  any vector $w\in\Sigma_{m}\Delta$ we have that
 $w\perp e_i e_j$ because  $\Delta\subset V_{e_i}\cap V_{e_j}$.
 For this vector
  we can construct
  a flat totally  geodesic triangle $e_ie_jx$ in $\mathcal A$ (as in Lemma~\ref{lem:mid})
  such that $\angle e_ie_jx=\pi/2$  and the vector $w$ is a
 tangent vector to this triangle at $m$. 
 By the definition of the map $\phi_i$ we have
 $d\phi_i(w)=\vr{e_jx}/|e_jx|$, let us note that $p(v)=\vr{e_jm}/|e_jm|$.
 We obtain that $d\phi_i(w)\perp p(v)$, then conclusion of 
 Lemma~\ref{lem:factor.sphere}  implies that the stabilizer
 $\Gamma_{e_j^0}$ can only fix or reflect  $e_i^0$.
 \qeds

 \section{Reflection property gives a lattice of singular points.}\label{secLattice}

In  this section we prove
Proposition~\ref{latticemain}.

For a point
$x\in\mathcal E$ we  denote the set of all adjacent
vertices by
$$\mathfrak S(x)=\{y\in\mathcal E|dim(V_x\cap V_y)=n-1\}.$$

The main technical point of the proof is the following:

\begin{thm}{Lemma}\label{4-th.point}
Let an action  $\Gamma\acts\R^n$ have the reflection
 property. Let  $x\in\mathcal E$,
$y,z\in\mathfrak S(x)$, and
 let us  denote $z^*=y+\overrightarrow{xz}$.
Then $$z^*\in\mathcal E.$$
\end{thm}
The proof of this lemma is in Section~\ref{pf4}.
In Section~\ref{prProp} we finish the proof of 
Proposition~\ref{latticemain}.

\subsection{Proof of Lemma~\ref{4-th.point}}\label{pf4}

We consider two cases.
First,  every element of the stabilizer $\Gamma_x$
  may reflect or fix  points $y$ and $z$ only simultaneously.
 The proof of this case is in Subsection~\ref{central} (see Lemma~\ref{lem:4th.point.general}(2)).
The other possibility  is  if there exists an element in $\Gamma_x$
 that reflects point $y$ and fixes point $z$,
 the proof for this case is in Subsection~\ref{noncentral}.
  
 For any two points $x,y\in\R^k$
we will denote by
$\mathfrak c_x: \R^k\to\R^k$ the
central symmetry with the center $x$
and by $\mathfrak c_{xy}: \R^k\to\R^k$
 the
 symmetry with an axis $xy$.
For points $x_1,\dots, x_l$ we
denote by $\<x_1,\dots, x_l\>$ the affine hull
of these points. Then in the first case 
$\Gamma_x|_{\<x,y,z\>}=\{id|_{\<x,y,z\>},\mathfrak c_x|_{\<x,y,z\>}\}$
and in the second case
$\Gamma_x|_{\<x,y,z\>}=\{id|_{\<x,y,z\>},\mathfrak c_x|_{\<x,y,z\>},
\mathfrak c_{xy}|_{\<x,y,z\>},\mathfrak c_{xz}|_{\<x,y,z\>}\}$.

\subsubsection{The  order of  $\Gamma_x|_{\<x,y,z\>}$ equals  2}\label{central}

For any two points $x,y\in \R^n$ we denote the stabilizer by
$\Gamma_{x,y}=\Gamma_x\cap\Gamma_y$ and by
$\Gamma^\#_{x,y}\subset O(n)$
the action of this stabilizer on the associate vector space $\R^n$. 
 For any tree points $x,y,z\in \R^n$ we denote the stabilizer of
 these points by
$\Gamma_{x,y,z}=\Gamma_x\cap\Gamma_y\cap\Gamma_z$.

First we prove  one auxiliary statement:
\begin{thm}{Lemma}\label{aux}
%$v\in \Sigma_{e_j}K_i$ so that $v\perp \uparrow_{e_j}^{e_i}$for every vector $v\in\Sigma_{M_{ij}}$ so hat$v\perp \uparrow_{M_{ij}}^{e_i}$the following are equivalent:1) there is an opposite v^*\in\Sigma_{M_{ij}}$: $|vv^*|=\pi$2) $\Gamma_{e_i^0e_j^0}(v)=v$
For a point $x\in\mathfrak S(y)$ and every vector $v\in\R^n$, if
$\Gamma^\#_{x,y}(v)=\{v\}$ then
 $\Gamma^\#_{y}(v)=\{v,-v\}$.

\end{thm}
\parit{Proof.} 
We can find points $x_1,\dots,x_{n-1}\in\mathfrak S(y)$,  such that vectors
$\overrightarrow{yx},\overrightarrow{yx_1},\dots,\overrightarrow{yx_{n-1}}$
are linearly independent.
Let $P=\{v\in\R^n|\Gamma_{x,y}^\#(v)=v\}$. Reordering if necessary we can 
assume that $x=x_0, x_1,\dots,x_k\in P$ and $x_{k+1},\dots,x_{n-1}\notin P$.
We know that 
$\Gamma^\#_y(\overrightarrow{yx_i})=\{\overrightarrow{yx_i},-\overrightarrow{yx_i}\}$
for every $i=0,\dots, n-1$.
Then considering the group action for the decomposition in our basis
$v=v^0\overrightarrow{yx_0}+\dots+v^{n-1}\overrightarrow{yx_{n-1}}$
we obtain that  for every $v\in P$
coordinates $v^{k+1}=\dots=v^{n-1}=0$, i.e.
$P=\<x_0,\dots,x_k\>$. Then for every $v\in P$ we have $\Gamma^\#_{y}(v)=\{v,-v\}$.
 \qeds

\begin{thm}{Lemma}\label{lem:4th.point.general}
In conditions of Lemma~\ref{4-th.point}
suppose additionally that $$\Gamma_x|_{\<x,y,z\>}=\{id|_{\<x,y,z\>},\mathfrak c_y|_{\<x,y,z\>}\}.$$

Then

1) $\Gamma_y|_{\<x,y,z\>}=\{id|_{\<x,y,z\>},\mathfrak c_y|_{\<x,y,z\>}\}$;

2) $z^*\in\mathcal E$.
 \end{thm}
\parit{Proof.}
1) Conditions of the lemma imply that
$\Gamma^\#_{x,y}(\{\overrightarrow{xz}\})=\{\overrightarrow{xz}\}$ and 
$ \overrightarrow{xz}=\overrightarrow{y{z^*}}$
then by Lemma~\ref{aux} $\Gamma_y(\overrightarrow{y{z^*}})=\{\overrightarrow{yz^*},-\overrightarrow{y{z^*}}\}$.
Then for every $\gamma\in\Gamma_y$ we have
 $\gamma(\overrightarrow{xz)}=\overrightarrow{xz}$ and
$\gamma(\overrightarrow{y{z^*}})=\overrightarrow{y{z^*}}$ or
 $\gamma(\overrightarrow{xz})=-\overrightarrow{xz}$ and
$\gamma(\overrightarrow{y{z^*}})=-\overrightarrow{y{z^*}}$.  Then (1) follows.

2) For every $v\neq 0$ we want to find
$\gamma\in\Gamma_{z^*}$ such that $\gamma(v)\neq v$. We consider two
possibilities.

First if $\Gamma_{x,y,z}(v)\neq\{v\}$ we 
can find the required element $\gamma\in\Gamma_{x,y,z}\subset\Gamma_{z^*}$.

If $\Gamma_{x,y,z}(v)=\{v\}$
then for arbitrary three elements $\gamma_x\in\Gamma_x\setminus\Gamma_{x,y,z}$, 
$\gamma_y\in\Gamma_y\setminus\Gamma_{x,y,z}$, $\gamma_z\in\Gamma_z\setminus\Gamma_{x,y,z}$ by Lemma~\ref{aux}
we will have $\gamma_x(v)=-v,\gamma_z(v)=-v,\gamma_z(v)=-v$.
Then $\gamma_x\circ\gamma_y\circ\gamma_z(v)=-v$ and
$\gamma_x\circ\gamma_y\circ\gamma_z(z^*)=z^*$.

\qeds

\subsubsection{The  order of  $\Gamma_x|_{\<x,y,z\>}$ equals  4}\label{noncentral}
In this subsection
 we are in conditions of  Lemma~\ref{4-th.point}
and assume that the order of the group
$\Gamma_x|_{\<x,y,z\>}$ is $4 $.

\begin{rdef}{Definitions}
Let vectors $v_1,v_2\in\R^n$
and $W$ be a subset of vectors in $\R^n$.
 The 
\emph{angle chain between $v_1$ and $v_2$ through the set $W$}
is the ordered set of vectors
$w_1,\dots,w_k\in W$ 
with the folowing property:
$$\angle(v_1, w_1)\neq\pi/2,
\angle(w_1, w_2)\neq\pi/2, \dots,\angle(w_{k-1}, w_k)\neq\pi/2, \angle(w_k, v_2)\neq\pi/2.$$

We say that two 
 $(n-1)$-faces of a convex polyhedron in $\R^n$
 are \emph{adjacent} if they have a common $(n-2)$-face.\end{rdef}

First 
we prove the following lemmas.

\begin{thm}{Lemma}\label{lem:G4}
Suppose we are
in condition of the Lemma~\ref{4-th.point}
and the order of the group
$\Gamma_x|_{\<x,y,z\>}$ is $4 $.
Then
there is
 no angle chain between
$\vr{xy}$ and $\vr{xz}$
through the set of vectors
$\{\vr{xt}|t\in\mathfrak S(x)\}$ 
and
$\vr{xy}\perp \vr{xz}$.

\end{thm}
\parit{Proof.}
It is sufficient to note that if there would be such an angle chain, then
the stabilizer $\Gamma_x$ could fix or reflect vectors $\vr{xy}$ and $\vr{xz}$
only simultaneously.   This contradicts to the fact that the order of the group
$\Gamma_x|_{\<x,y,z\>}$ is $4 $.
\qeds

\begin{rdef}{Notation}
We 
introduce the following notation
for  half spaces and hyperplanes determined by a vector $v\in\R^n$
or by an origin $x\in\R^n$ and a vector $v\in\R^n$:

$$H_{v}^{-} =\{w\in\R^n|\<w,v\>\le|v|^2\},\quad
H_{x,v}^{-} =\{y\in\R^n|\<\vr{xy},v\>\le|v|^2\} $$

$$H_{v}^{0} =\{w\in\R^n|\<w,v\>=|v|^2\} , \quad
H_{x.v}^{0} =\{y\in\R^n|\<\vr{xy},v\>=|v|^2\} .$$
\end{rdef}

We need the following observation in geometry of convex polyhedra.

\begin{thm}{Lemma}\label{lem:F}
Let an  $n$-dimensional
convex polyhedron
$F\subset\R^n$ be represented as an intersection of
half spaces
$F=\cap H_{v_i}^-$ for the set of vectors
$V=\{v_1,\dots,v_k\}\subset \R^n$
and suppose this set is minimal (or equivalently that
all intersections $F_i= H_{v_i}^{0}\cap F$ are hyperfaces in $F$ ).
Suppose in addition that
 there is
 no angle chain between
$v_1$ and $v_2$
through the set
$V$ and that
$v_1\perp v_2$.

Then faces
$F_1=H_{v_1}^{0}\cap F$ and $F_2=H_{v_2}^{0}\cap F$
are adjacent.% (i.e. $\dim(F_1\cap F_2)=n-2$).
\end{thm}

\parit{Proof.}

We consider a partition of $V$ into $3$ subsets.
$V_1$ is a subset
of those vectors in $V$ that have an angle chain to $v_1$
through $V$;
$V_2$ is a subset
of those vectors in $V$ that have an angle chain to $v_2$
through $V$
and $V_3=V\setminus(V_1\cup V_2)$.

Let $L=H_{v_1}^0\cap H_{v_2}^0$.
For any vector $w\in \R^n$ we denote by 
$w^L$ the orthogonal projection of $w$ to $L$.

Let $h_w=H^-_w\cap L$, then $h_w$ is a halfspace in $L$ with a normal vector $w^L$.
Set

$$I_1=\bigcap_{w\in V_1} h_{w},\quad
I_2=\bigcap_{w\in V_2} h_{w},\quad
I_3=\bigcap_{w\in V_3} h_{w}.$$
We claim that
\begin{enumerate}[(i)]
\item $\dim I_i=n-2$, for $i=1, 2, 3$;

\item for
$w_i\in V_i$  $w_j\in V_j$ we have
$w_i^L\perp w_j^L$ if $i\neq j$
\end{enumerate}
Let us first note,  that these two properties 
imply that $\dim (I_1\cap I_2\cap I_3)=n-2$.
This would imply the lemma, because
 $F_1\cap F_2=\cap_{w\in V}h_w=I_1\cap I_2\cap I_3$. 

Let us show (i) and (ii).

We define sets
$W_1=\{v\in V|\<v,v_2\>=0\}$ and
$W_2=\{v\in V|\<v,v_1\>=0\}$. 
Let note that: 
$V=W_1\cup W_2$,
$W_1\supset V_1$, 
$W_2\supset V_2$ and
$W_1\cap W_2\supset V_3$.

We define:
$$J_1=\bigcap_{w\in W_1}h_w,\quad J_2=\bigcap_{w\in W_2}h_w.$$
We  consider
$J_1$ as an intersection
$$J_1=(\bigcap_{w\in W_1}H_w^-\cap H_{v_1}^0)\cap L\subset H_{v_1}^0.$$
Let us note that:

1)$L$ is a hyperplane in $H_{v_1}^0$
with a normal vector $v_2$;

2) for every $w\in W_1$ the set $H_w\cap H_{v_1}^0$
is a half space in $H_{v_1}^0$ with a normal vector 
orthogonal to $v_2$

3)$\dim(\bigcap_{w\in W_1}H_w^-\cap H_{v_1}^0)=n-1$
(this follows from  the inclusion: $\bigcap_{w\in W_1}H_w^-\cap H_{v_1}^0\supset F_1$).

These three properties imply that $\dim J_1=n-1$.
By the same arguments $\dim J_2=n-1.$
We know that
$I_1\supset J_1$,  $I_2\supset J_2$ and
$I_3\supset J_1\cup J_2$, hence (i) follows. 

For any $w_1\in W_1$ and $w_2\in W_2$
the condition $w_1\perp w_2$ implies that $w_1^L\perp w_2^L$, hence (ii).

\qeds

For two points $x,y\in\mathcal E$
 we denote the common $(n-1)$-face of 
 polyhedra
 $V_x$ and $V_y$ by
 $V_{xy}=V_x\cap V_y$.  
 
 \begin{thm}{Lemma}\label{lem:4p}
Suppose that points
$x,y,z\in\mathcal E$,
$xy\perp xz$ and
 faces $V_{xy}$, $V_{xz}$ are adjacent.
 Suppose in addition that
 for any point 
$y^*\in\mathfrak S(y)\cap\<x,y,z\>$ we have
$yy^*\perp yx$.
Then the point $z^*=x+\vr{xy}+\vr{xz}\in\mathcal E$.

\end{thm}

\parit{Proof.}
Consider  the intersection $V_{xyz}=V_{xy}\cap V_{xz}$,
this intersection is an
$(n-2)$-face
of polyhedron $V_x$ and hence it is
an $(n-2)$-face
of polyhedron $V_y$ (because of the Voronoy decomposition structure). 
We can represent $V_y$
as an intersection:
$$V_{y}=\bigcap_{t\in\mathfrak S(y)}H_{y,\frac{1}{2}\vr{yt}}^{-}$$
Hence there 
 exist two points $y_1, y_2\in\mathfrak S(y)\cap\<x,y,z\>$
such that
$$V_{xyz}\subset H_{\frac{1}{2}\vr{yy_1}}^0\cap H_{\frac{1}{2}\vr{yy_2}}^0.$$
One of this points, say $y_1$ coincide with $x$. Then the other 
point  $y_2$ coincide with $z^*$.

\qeds

\parit{Proof of Lemma~\ref{4-th.point}}
We can represent $V_x$
as an intersection:
$$V_{x}=\bigcap_{t\in\mathfrak S(x)}H_{x,\frac{1}{2}\vr{xt}}^{-}$$
%=\cap_{Y\in\mathcal E\setminus X}H(X,Y)\eqlbl{V=H}.$$

Because of Lemma~\ref{lem:G4}
we can apply Lemma~\ref{lem:F}
to
$F=V_x$, 
$F_1=V_{xy}$,
$F_2=V_{xz}$ 
and obtain
that faces
$F_1=V_{xy}$,
$F_2=V_{xz}$ are adjacent.
Then we  apply Lemma~\ref{lem:4p}
 and
 Lemma~\ref{4-th.point} follows.
 \qeds

\subsection{}\label{prProp}
Here we finish the proof of
Proposition~\ref{latticemain}.
The proposition follows directly from
the next two lemmas.
\begin{thm}{Lemma}\label{latticebasis}
%For every $x\in\mathcal E$
The set $\mathcal E$ is a lattis with
 generating vectors $\{\overrightarrow{xy}\}_{y\in\mathfrak S(x)}$.

\end{thm}

\parit{Proof.}
We know that
$$V_{x}=\bigcap_{t\in\mathfrak S(x)}H_{x,\frac{1}{2}\vr{xt}}^{-}=
\bigcap_{t\in\mathcal E\setminus x}H_{x,\frac{1}{2}\vr{xt}}^{-}.\eqlbl{V=H}$$

%Let note that the reflection property implies that sets $V_X$,  $\mathfrak S(X)$ have central symmetry with center $X$.

 %can be reformulated as follows: for any $Y\in\mathfrak S(X)$ we have $\Gamma_X(XY)=\{XY,-XY\}$. 

First we show that
for every $x,y\in\mathcal E$

$$\mathfrak S(y)=\mathfrak S(x)+\overrightarrow{xy}.\eqlbl{S=S+xy}$$

Indeed Lemma~\ref{4-th.point} implies that
for any $x,y\in\mathcal E$ such that $y\in\mathfrak S(x)$
(further we call such points adjacent)
we have
$\mathfrak S(x)+\overrightarrow{xy}\subset\mathcal E$
and hence by \ref{V=H} 
$V_y\subset V_x+\overrightarrow{xy}$. 
Then changing $x$ and $y$
we obtain the equality \ref{S=S+xy} in this case.
For arbitrary  $x,y\in\mathcal E$ the equality \ref{S=S+xy} can be obtained
by joining $x$ and $y$ with a chain of adjacent points. 

To show Lemma~\ref{latticebasis} it is sufficient to prove that
for any $x,y,z\in\mathcal E$ we have
$x+\overrightarrow{xy}+\overrightarrow{xz}\in\mathcal E$
$x-\overrightarrow{xy}\in\mathcal E$.
The second inclusion needs 
 the central symmetry of the set $\mathfrak S(x) $, that follows from the reflection property.
 After this
 both inclusions
can be proved by using
\ref{S=S+xy} and joining correspondent points
 with a chain of adjacent points.
 \qeds

The next lemma
shows that we can reduce our generating set for
$\mathcal E$ to
the $n$ generating vectors with the same property that
$\Gamma$ reflects this vectors. 

\begin{thm}{Lemma}
Suppose 
$S\in\R^n$ is an $n$-lattis with a generating set $a_1,\dots, a_s$, 
$G$ is a subgroup of isometries of $\R^n$ 
which reflects this generating set.
Then there exists an $n$-generating set $b_1,\dots,b_n$ for $S$
such that  $G$
 reflects this generating set.

\end{thm}
\parit{Proof.}  It is known that there exists an $n$-generating set
for any $n$-lattice (sometimes it's called a short basis), so the problem is
to find  an $n$-generating set reflected by $G$.

We define an equivalence relation on the generating set $a_1,\dots,a_s$.
 We set $a_{i^*}\sim a_{j^*}$ if 
 there is an angle chain  connecting 
  $a_{i^*}$ and   $ a_{j^*}$
  through $\{a_1,\dots,a_s\}$.  
   %we can find sequence $a_{i^*}=a_{i_1},\dots,a_{i_k}=a_{j^*}$ so that $\angle (a_{i_m}, a_{i_{m+1}})\neq\pi/2$.
We denote by
$q_1,\dots,q_l$ the equivalence classes.
Let note that vectors from the different classes are mutually
orthogonal, equivalently vectors can be fixed or reflected by 
stabilizers $G_x$ (where $x\in S$) only simultaneously.

We fix a point $x\in S$.
For every equivalence class 
we consider  a lattice $$S_i=\{x+\sum_{a\in q_i}n_a a, n_a\in\mathbb Z\},$$
then
$$ S=\bigoplus_{i=1}^l\mathcal S^i.$$

 Let $d(i)$ be the dimension of the affine hull of $S_i$,
then we can choose  a $d(i)$-basis for $S_i$.  (We note that this basis
is independent on the choice of $x\in S$).
Vectors of this basis are reflected by $G$ because any stabilizer can act 
on every $S_i$ only identically or by central symmetry.
Then  the union of these $d(i)$-bases 
is an $n$-basis for $S$ reflected by $G$.

\qeds

\section{Comments and open questions.} 
Let 
$\Gamma\acts\RR^n$
be a discrete action on $\R^n$ by affine transformations.
Denote by $N(\Gamma)$ the number of orbits of isolated fixed points of
some subgroups in $\Gamma$.

Further, denote by $M(\Gamma)$ the number of maximal finite subgroups in $\Gamma$ up to conjugation.
Note that if $z$ is a singular point for the action $\Gamma\acts\R^n$ 
then the stabilizer of $z$ is a maximal finite subgroup of $\Gamma$.
It follows that $N(\Gamma)\le M(\Gamma)$.
Some maximal subgroups of $\Gamma$ might fix affine subspaces of positive dimension, therefore $M(\Gamma)$ might be strictly bigger that $N(\Gamma)$.
From Proposition~\ref{prop:affine+extreme} and Theorem~\ref{main2} % (\ref{main}), 
we have the following:

\begin{thm}{Corollary}
For any cocompact discrete action by affine transformations. $\Gamma\acts\RR^n$, we have $N(\Gamma)\le2^n$.
\end{thm}

We believe that the following stronger statement is true.

\begin{thm}{Conjecture}
For any cocompact discrete action by affine transformations $\Gamma\acts\RR^n$, we have $M(\Gamma)\le2^n$.
\end{thm}

There is a discussion on this conjecture, see \cite{mathoverflow}.

\end{document}